\documentclass[10pt]{article}
\usepackage{amsmath,enumerate,amsfonts,color,amssymb, verbatim, graphicx}

\graphicspath{ {./Figures/Dislocations/}{./Figures/QTensors-Blocks/}
{./Figures/NumericalMinimizationLDG_Energy/}
{./FiguresNumericalMinimizationHM_Energy/}}

\usepackage{ntheorem, float}
\usepackage{authblk}

\usepackage{tikz}

\def\ZZ{{\mathbb Z}}
\def\RR{{\mathbb R}}

\def\Sphere{{\mathbb S}}

\def\mcF{{\mycal F}}

\newtheorem{theorem} {\sc  Theorem\rm} [section]

\newtheorem{lemma} [theorem] {\sc  Lemma\rm}
\newtheorem{proposition} [theorem] {\sc  Proposition\rm}

\theorembodyfont{\upshape}
\newtheorem{remark}[theorem]{\sc  Remark\rm}

\def\bproof{\noindent{\bf Proof.\;}}
\def\eproof{\hfill$\square$\medskip}

\newcounter{marnote}

\DeclareFontFamily{OT1}{rsfs}{}
\DeclareFontShape{OT1}{rsfs}{m}{n}{ <-7> rsfs5 <7-10> rsfs7 <10-> rsfs10}{}
\DeclareMathAlphabet{\mycal}{OT1}{rsfs}{m}{n}

\def\tr{{\rm tr}}

\def\mcS{{\mycal{S}}}

\def\be{\begin{equation}}
\def\ee{\end{equation}}

\def\tr{{\rm tr}}

\def\mcS{{\mycal{S}}}

\def\mcE{{\mycal E}}

\newcommand{\defeq}{\stackrel{\scriptscriptstyle \text{def}}{=}}
\numberwithin{equation}{section}

\def\be{\begin{equation}}
\def\ee{\end{equation}}
\def\bea#1\eea{\begin{align}#1\end{align}}
\def\non{\nonumber}

\begin{document}
%\title[Point defects in liquid crystals]
 \title{Half-integer point defects in the $Q$-tensor theory of nematic liquid crystals  }
 \date{\today}
%\author{G. Di Fratta, JM Robbins, V. Slastikov, A. Zarnescu}

\author[1]{ G. Di Fratta}
\author[1]{JM Robbins}
\author[1]{V. Slastikov}
\author[2,3]{A. Zarnescu}
\affil[1]{School of Mathematics, University of Bristol, Bristol, UK}
\affil[2]{Department of Mathematics, University of Sussex, Falmer, UK}
\affil[3]{Institute of Mathematics ``Simion Stoilow", Bucharest, Romania}
%\affil[ ]{\textit {\{email1,email2,email3,email4\}@xyz.edu}}

%\address[G.\ Di Fratta]{University of Bristol, Bristol, UK}
%\email[Giovanni Di Fratta]{giovanni.difratta@bristol.ac.uk}
%\address[JM \ Robbins]{University of Bristol, Bristol, UK}
%\email[Giovanni Di Fratta]{J.Robbins@bristol.ac.uk}
%\address[V.\ Slastikov]{University of Bristol, Bristol, UK}
%\email[Valeriy Slastikov]{Valeriy.Slastikov@bristol.ac.uk}
%\address[A.\ Zarnescu]{University of Sussex, Falmer, UK}
%\email[Arghir Zarnescu]{A.Zarnescu@sussex.ac.uk}

\maketitle

\begin{abstract} 
We investigate prototypical  profiles of point defects in two dimensional liquid crystals within the framework of  Landau-de Gennes theory. Using %suitable, defect-type boundary conditions
boundary conditions characteristic of %director-field 
defects of index $k/2$, we  find a  %class of  
critical point of the  Landau-de Gennes energy that is characterised by a system of ordinary differential equations. In the deep nematic regime, $b^2$ small, we prove that this critical point is the unique  global minimiser of the Landau-de Gennes energy. {\color{black} For the case $b^2=0$,} we investigate in greater detail the regime of vanishing elastic constant $L \to 0$, where we obtain three explicit point defect profiles, including the global minimiser.
\end{abstract}

% \maketitle
 
 \section{Introduction}
 
 Defect structures are among the most important and visually striking  patterns associated with nematic liquid crystals. These  are observed when passing polarised light through a liquid crystal sample and are characterised by sudden, localised changes in the intensity and/or polarisation of the light \cite{chandra,degennes}. Understanding the mechanism that generates defects and predicting their appearance and stability is one of the central objectives of any liquid crystal theory. %The existing major continuum theories of nematics, namely Oseen-Frank, Ericksen and de Gennes', do not seem to have a common interpretation of the defects \cite{ericksen, dgdefects, virga}.
 % What is nevertheless understood experimentally is that the defects can be classified depending on their dimension as either point, line  or wall  defects.
 
 %JR.  1) We do not discuss Ericksen, so would not mention it explicitly.  2) As the models are different, you would expect the interpretation of defects to be different, right?  3) I've added a reference to Kleman's book. 
 
\begin{figure}[ht]
\begin{minipage}[b]{0.45\linewidth}
\centering
\includegraphics[width=0.8\textwidth]{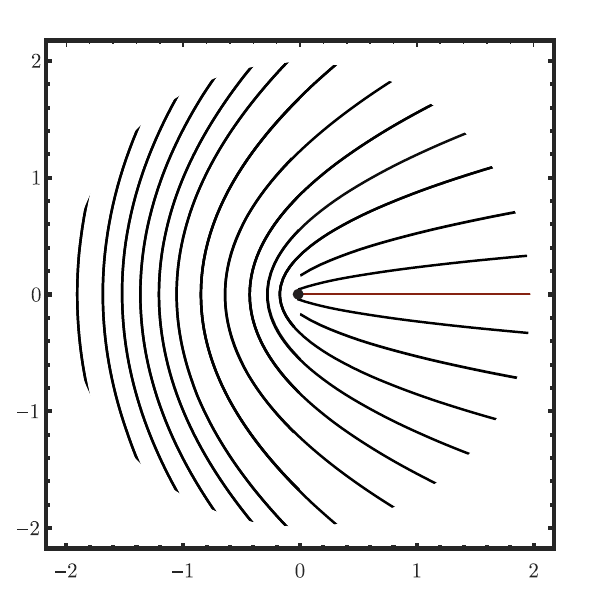}
%\caption{default}
\end{minipage}
\hspace{0.5cm}
\begin{minipage}[b]{0.45\linewidth}
\centering
\includegraphics[width=0.765\textwidth]{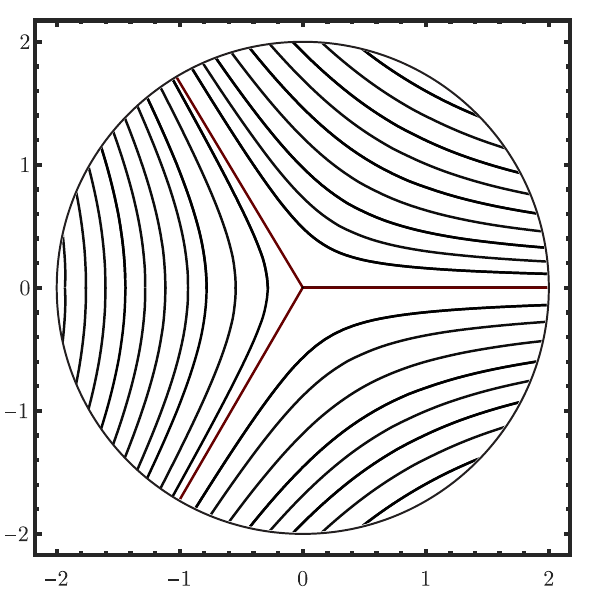}
%\caption{default}
\end{minipage}
\caption{Defects of index $\frac{1}{2}$ (left) and $-\frac{1}{2}$ (right) }
\label{fig:fig1}
\end{figure}

\begin{figure}[ht]
\begin{minipage}[b]{0.45\linewidth}
\centering
\includegraphics[width=0.8\textwidth]{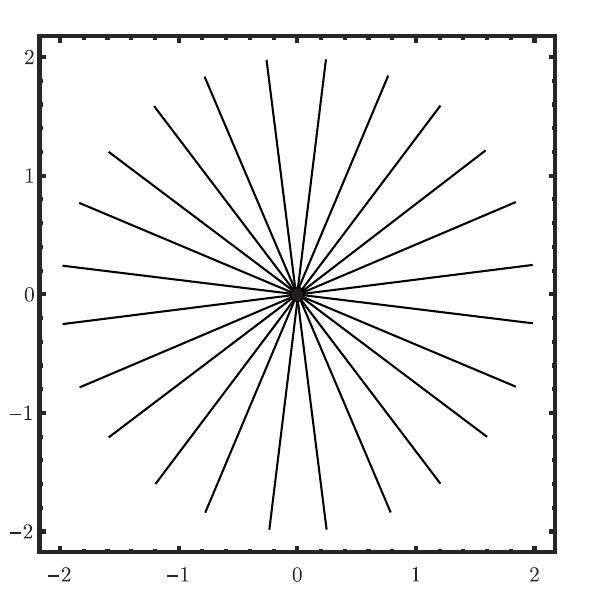}
%\caption{default}
\end{minipage}
\hspace{0.5cm}
\begin{minipage}[b]{0.45\linewidth}
\centering
\includegraphics[width=0.81\textwidth]{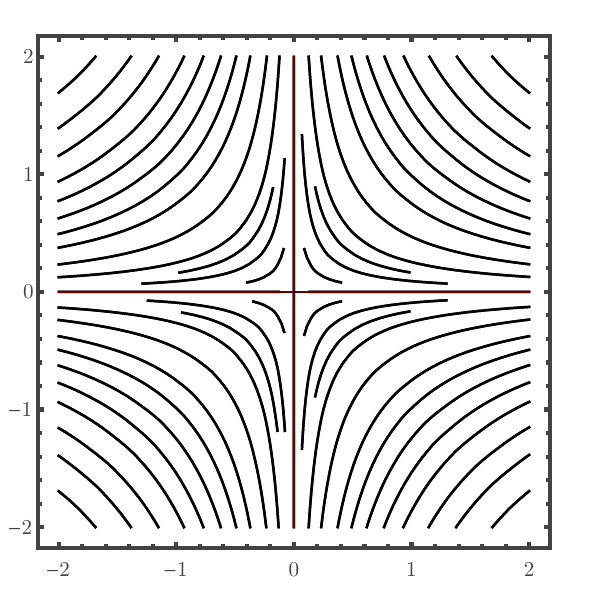}
%\caption{default}
\end{minipage}
\caption{Defects of index ${1}$ (left) and $-{1}$ (right) }
\label{fig:fig2}
\end{figure}

\begin{figure}[ht]
\begin{minipage}[b]{0.45\linewidth}
\centering
\includegraphics[width=0.8\textwidth]{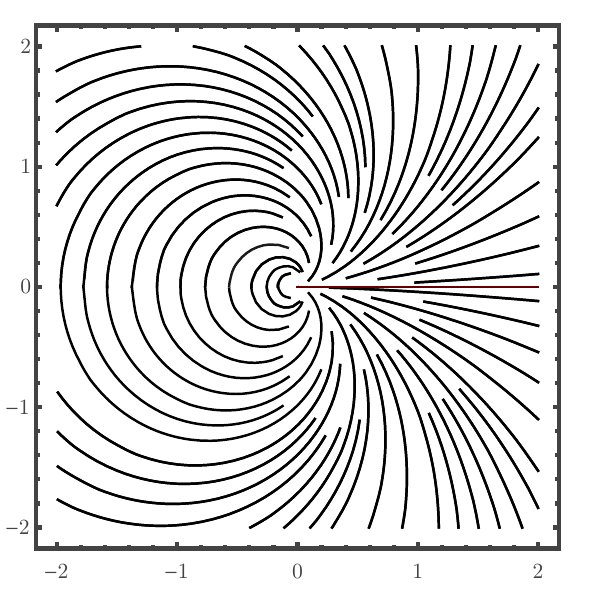}
%\caption{default}
\end{minipage}
\hspace{0.5cm}
\begin{minipage}[b]{0.45\linewidth}
\centering
\includegraphics[width=0.81\textwidth]{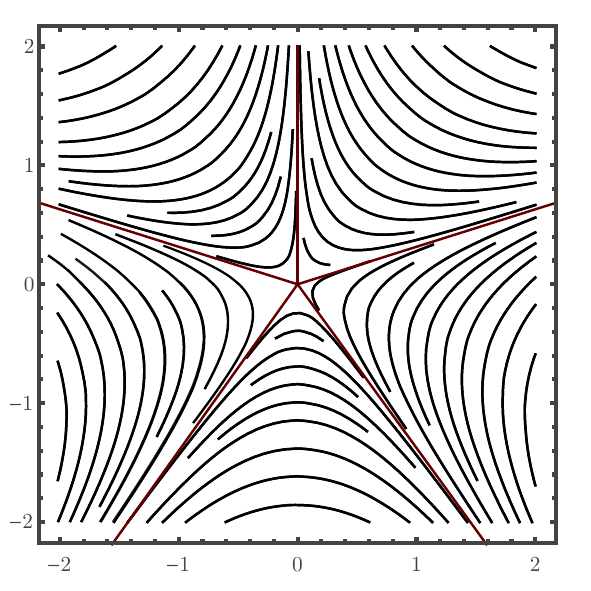}
%\caption{default}
\end{minipage}
\caption{Defects of index $\frac{3}{2}$ (left) and $-\frac{3}{2}$ (right) }
\label{fig:fig3}
\end{figure}

\begin{figure}[ht]
\begin{minipage}[b]{0.45\linewidth}
\centering
\includegraphics[width=0.8\textwidth]{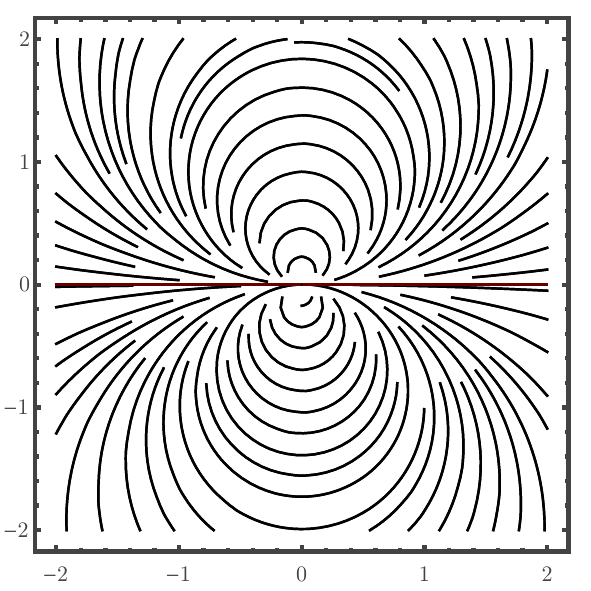}
%\caption{default}
\end{minipage}
\hspace{0.5cm}
\begin{minipage}[b]{0.45\linewidth}
\centering
\includegraphics[width=0.81\textwidth]{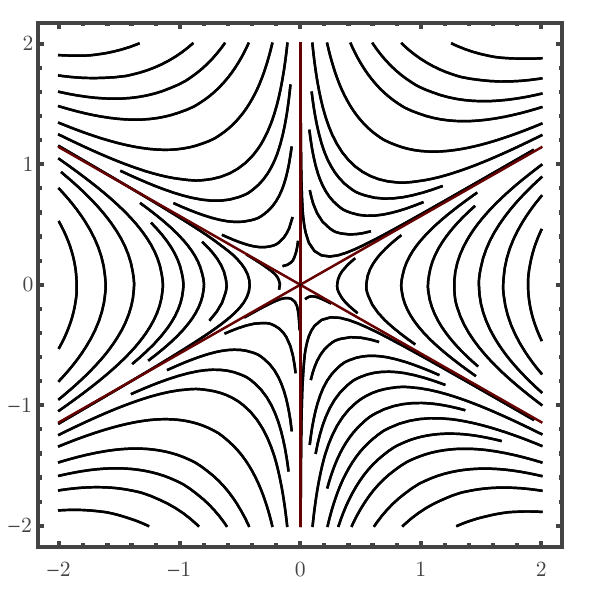}
%\caption{default}
\end{minipage}
\caption{Defects of index ${2}$ (left) and $-{2}$ (right)  }
\label{fig:fig4}
\end{figure}

 The mathematical 
 characterisation of defects depends  on the underlying model \cite{ericksen, degennes, kleman,  virga}.
In the Oseen-Frank theory, nematic liquid crystals are described by a vector field ${\bf n}$ defined on a domain $\Omega\subset\RR^d$   taking values in $\Sphere^{d-1}$ ($d=2,3$), which describes the mean local orientation of the constituent particles. 
%
%JR. Dimension of target space isn't always one less than dimension of domain; it's possible to have 2d domain (thin film) with n taking values in S^2.  But it's simpler to leave the text as is.
%
%According to this theory 
Defects correspond to discontinuities  %the vector field 
in ${\bf n}$ \cite{chandra, klemanlavrentovich, virga} and may be classified topologically. %. In order to classify the defects one assigns a topological invariant, or charge to every defect.  
For example, for planar vector fields %$\bf n$ 
in two-dimensional domains (i.e., $d=2$ above), point defects may be characterised by %an integer that specifies 
the number of times $\bf n$ rotates through $2\pi$ as an oriented circuit around the defect is traversed.
%or instance, it is possible to characterise the strength (or the charge) of the point defect by a special number  $k \in \ZZ$. This number corresponds to how many times director ${\bf n}$ rotates by $2\pi$ when one goes around %defect core, and the sign of $k$ indicates a direction of rotation \cite{chandra, krlav, toukle}. 
%
% JR.  Have removed k from the preceding.  Previously, k meant two different things here and in the following.
%
%In Ericksen's  model liquid crystal is described using the order parameter $s {\bf n}$, with ${\bf n}$ defined as in Oseen-Frank theory and a function $s\, : \, \Omega \to \RR$. Defects in this theory correspond to zeroes of $s$ and their strength is defined exactly as in the Oseen-Frank theory. 
For nonpolar nematic liquid crystals,  ${\bf n}$ and  $-{\bf n}$ are physically equivalent; in this case it is more appropriate %therefore more appropriate 
to regard ${\bf n}$ as taking values in $\RR \mathbb P^{d-1}$ rather than $\Sphere^{d-1}$. 
The classification of point defects in two dimensions then allows for  both integer and half-integer indices $k/2$, %$\frac{k}{2}$, 
$k \in \ZZ$ \cite{ballzarnescu, chandra, klemanlavrentovich}, as ${\bf n}$ is constrained to turn through a multiple of $\pi$ rather than $2\pi$ on traversing a circuit. Prototypical examples of such defects are shown in Figures~\ref{fig:fig1} -- \ref{fig:fig4}.

%Major 
A deficiency of the Oseen-Frank theory is that point defects in two dimensions, which are observed experimentally, %  defects %usually 
are predicted to have infinite energy; moreover, the theory does not allow for half-integer indices (see \cite{ballzarnescu, degennes}).  These shortcomings are addressed %One possibility to obtain qualitatively similar patterns for half-integer defects, as observed experimentally, is to use 
by the more comprehensive Landau-de Gennes $Q$-tensor theory \cite{degennes}. In this theory, the order parameter describing the liquid crystal system takes values in the space of $Q$-tensors (or $3\times 3$ traceless symmetric matrices),
$$
\mcS_0\defeq\{ Q\in \RR^{3\times 3},\, Q=Q^t,\,\tr (Q)=0\}.
$$ 
%According to Landau-de Gennes theory 
Equilibrium configurations of liquid crystals correspond to local minimisers of the Landau-de Gennes energy, which in its simplest form is given by% that can be written (in the simplest form) as
\be\label{ener}
 \mathcal{F}[Q]\defeq\int_{\Omega} \left\{ \frac{L}{2}|\nabla Q(x)|^2-\frac{a^2}{2}\tr(Q^2)-\frac{b^2}{3}\tr(Q^3)+\frac{c^2}{4}\left(\tr(Q^2)\right)^2 \right\} \,dx.
\ee
Here $Q \in \mcS_0$, $L>0$ is the elastic constant, and $a^2, c^2 > 0$,  $b^2 \geq 0$ are material parameters which may depend on temperature %depending on a physical regime 
(for more details see \cite{degennes}).

\begin{figure}[ht]
\begin{minipage}[b]{0.45\linewidth}
\centering
\includegraphics[width=1\textwidth]{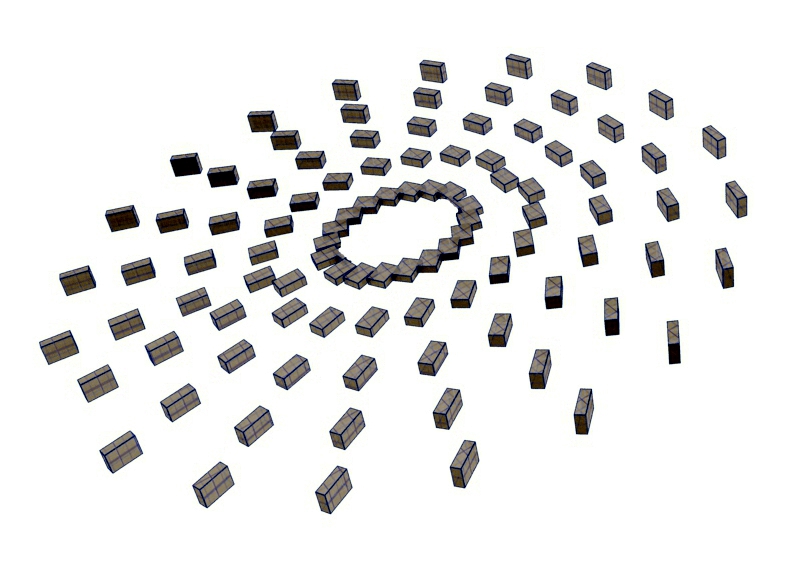}
%\caption{default}
\end{minipage}
\hspace{0.5cm}
\begin{minipage}[b]{0.45\linewidth}
\centering
\includegraphics[width=1.01\textwidth]{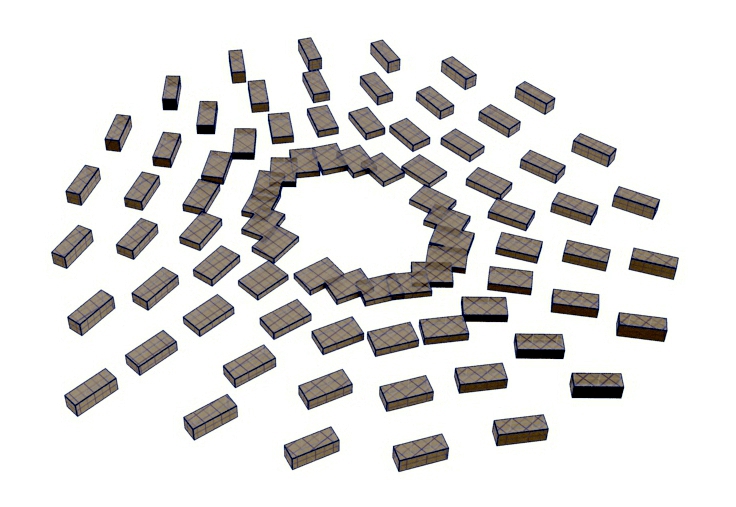}
%\caption{default}
\end{minipage}
\caption{$Q$-tensor defect of index $\frac{1}{2}$ (left) and $-\frac{1}{2}$  (right)  }
\label{fig:fig5}
\end{figure}

\begin{figure}[ht]
\begin{minipage}[b]{0.45\linewidth}
\centering
\includegraphics[width=1\textwidth]{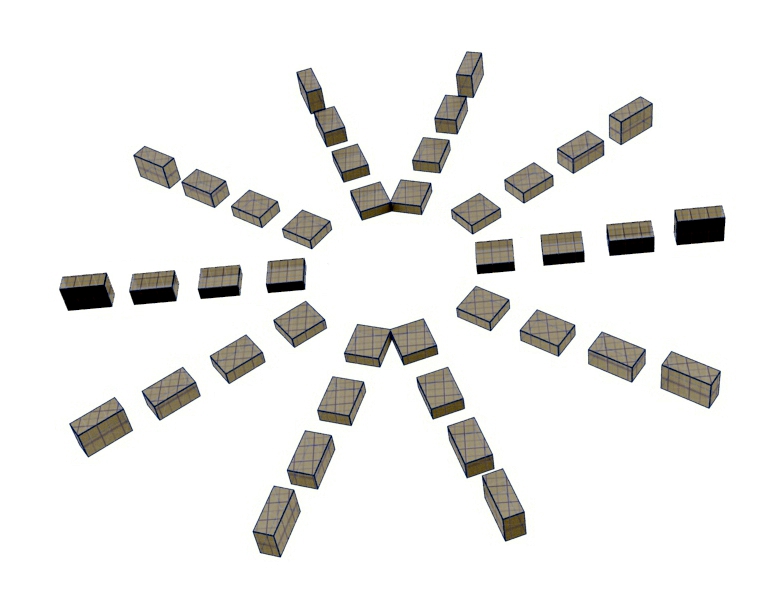}
%\caption{default}
\end{minipage}
\hspace{0.5cm}
\begin{minipage}[b]{0.45\linewidth}
\centering
\includegraphics[width=1.01\textwidth]{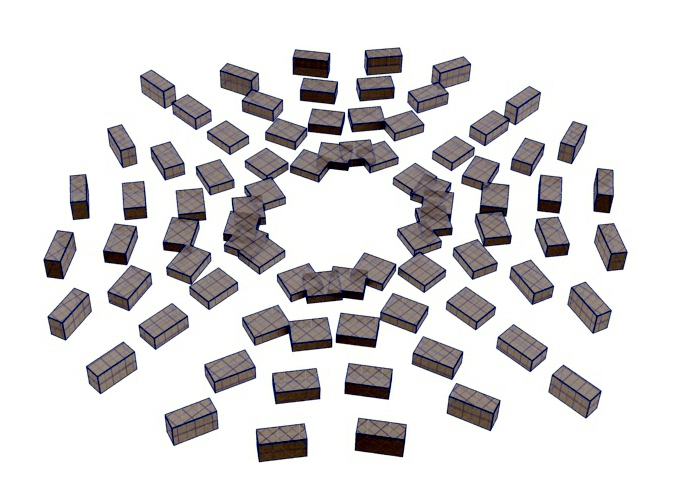}
%\caption{default}
\end{minipage}
\caption{$Q$-tensor defect of index ${1}$ (left) and $-{1}$ (right)  }
\label{fig:fig6}
\end{figure}

%%%%CHANGE
One can visualise $Q$-tensors as {\color{black} parallelepipeds} whose axes are parallel to the eigenvectors of  $Q(x)$ %the tensor 
with lengths given by the eigenvalues \cite{chopar}.\footnote{The careful reader will note that $\tr (Q) = 0$ implies that the eigenvalues cannot all be positive. %, since $\tr Q the sum of the eigenvalues is zero. 
In order to obtain positive lengths for the axes, %-tensors with bounded eigenvalues,
 %it suffices to 
 we add to each eigenvalue a sufficiently large positive constant (we assume the eigenvalues of  $Q$ are bounded).} % in order to obtain positive lengths for the axes the parallelepipeds}
Figure~\ref{fig:fig5}  displays defects of  index $\pm\frac12$ using this representation, and Figure~\ref{fig:fig6} displays defects of  index $\pm 1$.\footnote{The figures represent the numerically computed solutions of \eqref{ODEsystem}, \eqref{bdrycond} for $k=\pm 1,\pm 2$.}

\par The goal of this paper is a rigorous study %prototypical 
 of point defects in liquid crystals in two-dimensional domains using Landau-de Gennes theory.  We investigate equilibrium configurations in the disk $\Omega=\{ (x,y)\, :\, x^2 + y^2 < R\} $ subject to boundary conditions characteristic of prototypical defects, namely that on $\partial \Omega = \{(R\cos\phi, R\sin\phi)\}$, $Q$ is proportional to
\[ Q_k = \left( n \otimes n -\frac{1}{3} I \right), \ n =(\cos ({\textstyle\frac{k}{2}} \phi) , \sin ({\textstyle\frac{k}{2}} \phi) , 0).\]

%chosen to correspond to the physically observed defects of strength $\frac{k}{2}$, $k \in \ZZ$ (see Figures~\ref{fig:fig1} -- \ref{fig:fig4}).  
%In order to do this 
We first introduce an %following 
ansatz %in polar coordinates
\be\label{anY}
Y= u(r) \sqrt{2} \left (n(\varphi)\otimes n(\varphi)-\frac{1}{2}I_2\right) + v(r)  \sqrt{\frac{3}{2}} \left( e_3\otimes e_3-\frac{1}{3}I \right),
\ee
and note that $Y$ satisfies the Euler-Lagrange equations \eqref{eq:EL} for % of the critical points of 
the Landau-de Gennes energy \eqref{ener} provided that $(u,v)$ satisfies a system of ODEs given by \eqref{ODEsystem}, \eqref{bdrycond}. It follows that for all parameters $L, a, b, c$,  the ansatz 
$Y$ %defined in \eqref{anY} 
is a critical point of the energy.% \eqref{ener}.

Next, we show that for every $k \in \ZZ$, the critical point $Y$ is actually the unique global minimiser of the energy \eqref{ener} in the low-temperature  regime, i.e.~for $b^2$ sufficiently small. % ($b^2$ is small enough). 
{Equivalently, in this regime, $Y$ describes the unique ground state configuration for a two-dimensional index-$k$ point defect.}
% For this reason, we regard $Y$ as the {\it universal profile} for point defects in Landau-de Gennes theory in two dimensions in this regime.  %of the liquid crystal system. 
{\color{black} In general, it is very difficult to find a global minimizer of a non-convex energy. In this case we can deal with the nonlinearity using properties of the defect profile $(u,v)$ and the Hardy decomposition trick \cite{INSZ1}. 
Similar ideas to prove global minimality are used in \cite{shirokoffchoksi} for a problem in diblock copolymers. 
}

{\color{black} In the case $b^2=0$,} we also study the regime of vanishing elastic constant $L \to 0$ (see the appendix of \cite{nz} for a discussion of the  physical relevance of this regime) and show that it leads to a harmonic map problem for $Y$. We find three explicit solutions  -- two biaxial and, for even $k$, one uniaxial -- and show that one of the biaxial solutions is the unique global minimiser of \eqref{ener}.  The uniaxial critical point is  analogous to the celebrated ``escape in third dimension" solution of Cladis and Kl\'eman \cite{coaxialbbh, cladiskle}. 

%The problem of the liquid crystal defect 
The profile and stability of liquid crystal defects have been extensively studied in the mathematics literature \cite{danpatty, coaxialbbh, biscarivirga, canavaro, fatkul, golovaty, duvan_apala,   INSZ2, INSZ3, INSZ1, coaxialvirga, MakGar1}. Let us briefly mention  a few papers which bear directly on the present work. 
%One of the first results in this direction is the paper by Mironescu \cite{mironescu}. He investigated the stability of vortex-type defect in two-dimensional balls of finite radius depending on the winding number $k$. The main result is the stability of the defect for $k=1$ and instability for $k \geq 2$.ss
%{\Large A couple of sentences on Alama, Bronsard, Mironescu paper}
%\par 
In  \cite{coaxialvirga} the problem  of investigating equilibria of liquid crystal systems in cylindrical domains (effectively 2D disks) was studied numerically for the Landau-de Gennes model under homeotropic boundary conditions (i.e., $k=2$ above),  subject to the so-called Lyuksutov constraint $\textrm{tr}(Q^2)=a^2/c^2$. % (for $L>0$) constraint referred to in the physical literature as the  ``Lyuksutov's constraint''. 
The authors compare three different solutions of this model corresponding to ``planar positive", ``planar negative" and ``escape in third dimension". They numerically explore the energies of these solutions and find a crossover between the ``planar negative" and ``escape in third dimension" solutions depending on the parameters $b$ and $L$. For $b=0$, the  ``planar negative" solution is found to have lower energy than the other two.
%In the context of the Oseen-Frank theory the problem of stability of nematic states between two coaxial cylinders has been studied in \cite{coaxialbbh} and there it was related to the celebrated phenomenon of ``escape into the third dimension''. Its analogue in the $Q$-tensor framework was considered in \cite{biscarivirga}. 

In recent papers \cite{INSZ2, INSZ3, INSZ1} the radially symmetric 3D point defect, the so-called melting hedgehog, was studied  within the framework of Landau-de Gennes theory. The authors investigate the profile and stability of the defect as a function of the material constants $a^2, b^2, c^2$. In particular, it is shown that for $a^2$ small enough the melting hedgehog is locally stable, while for $b^2$ small enough it is unstable. We utilise some ideas introduced in the liquid crystal context in these papers to derive our present results.

The paper is organised as follows: The  mathematical formulation of the problem is given in section 2. %using $Q$-tensor theory. 
In section 3 we introduce an ansatz $Y$ satisfying  boundary conditions 
 %we obtain a critical point of the Landau-de Gennes energy that describes the profile of the
  characteristic of a point defect of index $k/2$, and show that Euler-Lagrange equations simplify from a system of PDEs to a system of two ODEs. We establish the existence of a solution of this system of ODEs, and thereby prove the existence of a critical point of the Landau-de Gennes energy.

In section 4 we study qualitative properties of the solution in the infinitely-low temperature regime, i.e.~for $b^2=0$. We study separately the case of fixed $L>0$ and the limit  $L \to 0$. The main result for fixed $L$  is that for all $k\in \ZZ$, $Y$ is the unique global minimiser of the Landau-de Gennes energy 
%  global minimality of the %special
%solution $Y$ with respect to variations in 
over $H^1(\Omega, \mcS_0)$.  Thus, for $b^2$ sufficiently small, %, and indeed for sufficiently small $Bb 
$Y$ describes the unique ground state for point defects in 2D liquid crystals. %This holds for all $k \in \ZZ$.  In particular, we show that for $b^2=0$,  defects of all strengths $\frac{k}{2}$, $k \in \ZZ$ are global minimizers of the Landau-de Gennes energy and therefore correspond to the universal profiles of 2D point defects of liquid crystal system. 
In the limit $L \to 0$, we derive the corresponding harmonic map problem and explicitly find three solutions -- two biaxial and, for even $k$, one uniaxial. We show that one of the biaxial solutions, $Y_-$, is the unique global minimiser of the Dirichlet energy.  Section 5 contains a discussion of the results and an outlook on further work.
% Section 5 is devoted to the study of the stability of the special solutions when $b^2$ is small enough. We show that  defects of all strengths $\frac{k}{2}$, $k \in \ZZ$ are local minimizers of the Landau-de Gennes energy and therefore correspond to stable equilibrium patterns of liquid crystal system. 

\section{Mathematical formulation of the problem}

%We are interested in studying equilibriums of liquid crystal systems on a two-dimensional disc of radius $R$.

We consider the following Landau-de Gennes energy functional on a two-dimensional domain $\Omega \subset \RR^2$, 
\be \label{FQ} 
 \mathcal{F}[Q]\defeq\int_{\Omega} \frac{1}{2}|\nabla Q(x)|^2 +\frac{1}{L} f(Q)\,dx, \quad Q\in H^1(\Omega; \mcS_0).  % \ \Omega \subset \RR^2.
\ee 
Here $L>0$ is a positive elastic constant, $\mcS_0$ denotes the set of $Q$-tensors defined by
$$\mcS_0\defeq\{ Q\in \RR^{3\times 3},\, Q=Q^t,\,\tr (Q)=0\}$$ %is the set of $Q$-tensors. T
and the bulk energy density $f(Q)$ is given by %has the following form
$$
f(Q)= -\frac{a^2}{2}|Q|^2-\frac{b^2}{3}\tr(Q^3)+\frac{c^2}{4}|Q|^4,
$$
where $a^2, c^2 >0$ and  $b^2 \geq 0$ are material parameters %constants 
and $|Q|^2 \defeq\textrm{tr} (Q^2)$.

We are interested in studying  critical points and local minimisers of the energy \eqref{FQ} for $\Omega = B_R$,  where $B_R\subset \RR^2$ is
%(here $B_R \subset \RR^2$ denotes 
the disk of radius $R < \infty$ centered at $0$, such that $Q$ satisfies boundary conditions %has suitable boundary conditions 
%on $\partial B_R$ 
corresponding to a point defect at $0$ of index $k/2$.  Specifically, we define  %we require that
% defined as follows  On $\partial B_R$, the director field $n$  % = \{ (R\cos\phi, R\sin\phi)$\} 
%for a point defect of index $k/2$, where $k \in \ZZ$, is given in polar coordinates   by
%\be
%\label{eq: n(phi)}
%n(\varphi)=\left(\cos\left(\frac{k}{2}\varphi\right),\sin\left(\frac{k}{2}\varphi\right),0\right).% \quad k\in\mathbb{Z} \setminus \{ 0 \},
%\ee
%The corresponding $Q$-tensor is given by
%$n =\left( \cos\left(\frac{k}{2} \varphi \right),   \sin\left(\frac{k}{2} \varphi \right), 0 \right)$. The corresponding $Q$-tensor is given by
%, $k\in \ZZ$, gives rise to the $Q$-tensor field  of strength $k/2$ for $k\in \Zz$. typical defect of degree $\frac{k}{2}$ ($k \in \ZZ$, $k\neq 0$) has the following form
\be\label{def:qkbdry}
Q_k(\varphi) =  \left( n(\varphi) \otimes n(\varphi) - \frac{1}{3} I \right),
\ee
where
\be \label{eq: n(phi)}
n(\varphi)=\left(\cos\left(\frac{k}{2}\varphi\right),\sin\left(\frac{k}{2}\varphi\right),0\right), \quad k\in\mathbb{Z} \setminus \{ 0 \},
\ee
%$n =\left( \cos\left(\frac{k}{2} \varphi \right),   \sin\left(\frac{k}{2} \varphi \right), 0 \right)$ is a unit vector expressed in terms of the polar coordinate $\phi$ on $\partial B_R$,  
and $I$ is the $3\times 3$ identity matrix.
%where $n =\left( \cos\left(\frac{k}{2} \varphi \right),   \sin\left(\frac{k}{2} \varphi \right), 0 \right)$ is a unit vector expressed in the polar coordinates and 
%where $I$ is the $3\times 3$ identity matrix.  
The boundary condition  is then taken to be %We require that %impose the following boundary conditions on $Q$:
\be \label{BC}
Q(x) = s_+Q_k(\varphi) \quad \hbox{for all } x \in \partial B_R, 
\ee
where $x = (R\cos\phi, {\color{black} R\sin \phi})$ and 
\be \label{def:s_+}
s_+ = \frac{ b^2 + \sqrt{b^4+24 a^2 c^2}}{4 c^2}.
\ee The value of $s_+$ is chosen  so that $s_+ Q_k$ minimizes $f(Q)$. 
Critical points of the energy functional satisfy  the Euler-Lagrange equation:
\be\label{eq:EL}
L \Delta Q=-a^2 Q-b^2[Q^2-\frac 13|Q|^2I]+c^2Q|Q|^2 \hbox{ in } B_R, \ \  Q=s_+Q_k \hbox{ on } \partial B_R,
\ee 
where the term $b^2\frac 13|Q|^2 I$ accounts for the constraint $\tr(Q)=0$.

\section{Existence of special solutions}

In general it is difficult to find %solutions of the Euler-Lagrange equation for the
critical points of the Landau-de Gennes energy.  However, due to symmetry we are able to find a special class of solutions of  the Euler--Lagrange equation \eqref{eq:EL}. % that turn out to describe the defect of degree $\frac{k}{2}$, for any ($k \in \ZZ$, $k \neq 0$)

We consider the following ansatz, expressed in polar coordinates $(r,\varphi)\in (0, R) \times [0,2\pi]$:
\be\label{ansatz}
Y(r,\varphi)=u(r)F_n (\varphi)+v(r)F_3,
\ee 
where %we denote
\be
F_n(\varphi)\defeq \sqrt{2} \left (n(\varphi)\otimes n(\varphi)-\frac{1}{2}I_2\right), \quad
F_3\defeq \sqrt{\frac{3}{2}} \left( e_3\otimes e_3-\frac{1}{3}I \right),
\ee 
%where
%\be
%n(\varphi)=\left(\cos\left(\frac{k}{2}\varphi\right),\sin\left(\frac{k}{2}\varphi\right),0\right), \quad k\in\mathbb{Z} \setminus \{ 0 \},
%\ee
$n(\varphi)$ is given by \eqref{eq: n(phi)} and $I_2 = e_1\otimes e_1 + e_2 \otimes e_2$ ($e_i$ denotes the standard basis vectors in $\RR^3$).
It is  straightforward %calculation 
to check that 
 $|F_n|^2=|F_3|^2=1$ and $\tr(F_nF_3)=0$, so that $Q_k$ may be expressed as
%t follows that
%
 %and for $Q_k$ as in \eqref{def:qkbdry} we have
$$
Q_k (\varphi) =  \frac{1}{\sqrt{2}} F_n(\varphi) - \frac{1}{\sqrt{6}} F_3 .
$$
It follows that  $Y(r,\phi)$ % we obtain that 
%$$Y(R,\varphi) = s_+ Q_k(\varphi)$$ 
satisfies the boundary conditions \eqref{BC} provided
\be \label{eq: bc at R}
u(R) =\frac{1}{\sqrt{2}}s_+, \quad v(R)=-\frac{1}{\sqrt{6}} s_+.
\ee 
%$Y(R,\phi)$ % we obtain that 
%$$Y(R,\varphi) = s_+ Q_k(\varphi)$$ 
%satisfies boundary conditions \eqref{BC}.

\begin{remark}
For $k=2$, $Y(r, \varphi)$ satisfies hedgehog boundary conditions  (see Figure~\ref{fig:fig6}, left), while for $k=\pm1$, $Y$ satisfies boundary conditions corresponding to a  defect of index $\pm \frac12 $ \cite{chandra, klemanlavrentovich}.
%with topological charge plus or minus one. 
The  $-\frac12 $-defect is also called a $Y$-defect because of its shape (see Figure~\ref{fig:fig5}, right).

%\begin{figure}[ht]
%\begin{minipage}[b]{0.45\linewidth}
%\centering
%\includegraphics[width=1\textwidth]{Newk12.jpg}
%%\caption{default}
%\end{minipage}
%\hspace{0.5cm}
%\begin{minipage}[b]{0.45\linewidth}
%\centering
%\includegraphics[width=1.01\textwidth]{Newkm12.jpg}
%%\caption{default}
%\end{minipage}
%\caption{$Q$-tensor defect of strength $\frac{1}{2}$ (left) and $-\frac{1}{2}$  (right)  }
%\label{fig:fig6}
%\end{figure}
\end{remark}

\par We would like to show that the ansatz \eqref{ansatz} satisfies the Euler-Lagrange equation \eqref{eq:EL} provided $u(r)$ and $v(r)$ satisfy a certain system of ODEs. It is straightforward to check that
\bea\label{deltatopolar}
\Delta Y=\left(u''(r)+\frac{u'(r)}{r}-\frac{k^2u(r)}{r^2}\right)F_n(\varphi)+\left(v''(r)+\frac{v'(r)}{r}\right)F_3
\eea
% On the other hand we have:
and
\bea\label{q2ansatz}
Y^2= -\sqrt{\frac{2}{3}} uv F_n (\varphi) + \frac{1}{\sqrt{6}} \left( - u^2+ v^2 \right) F_3 +\frac{1}{3} |Y|^2 I, \quad
 |Y|^2=u^2+v^2.
\eea
Substituting (\ref{ansatz}), (\ref{deltatopolar}) and (\ref{q2ansatz}) into (\ref{eq:EL}) we obtain
\bea
&\left(u''(r)+\frac{u'(r)}{r}-\frac{k^2u(r)}{r^2}\right)F_n(\varphi)+\left(v''(r)+\frac{v'(r)}{r}\right)F_3\non\\
&= \frac{1}{L}\left(-a^2u+\sqrt{\frac{2}{3}}b^2uv+{c^2u}\left( u^2+ v^2\right)\right){\color{black} F_n(\varphi)}\nonumber\\
&+\frac{1}{L}\left(-a^2v-\frac{1}{\sqrt{6}} b^2\left(  -u^2+ v^2 \right)+{c^2v}\left( u^2+ v^2\right) \right)F_3.
\eea
Taking into account that the matrices $F_n(\varphi),F_3$ are linearly independent for any $\varphi\in [0,2\pi]$ we obtain 
the following coupled system of ODEs for $u(r)$ and $v(r)$: % that provides a solution to \eqref{eq:EL}
%\be\label{ODEsystem}
%\left\{\begin{array}{l}
%u''+\frac{u'}{r}-\frac{k^2u}{r^2}=\frac{u}{L}\left[-a^2+\sqrt{\frac{2}{3}} b^2 v+c^2\left( u^2+ v^2\right)\right]\\
%v''+\frac{v'}{r}=\frac{v}{L}\left[-a^2-\frac{1}{\sqrt{6}}b^2 v+c^2\left( u^2+ v^2\right) \right] + \frac{1}{L \sqrt{6}} b^2 u^2  
%\end{array}\right. \hbox{ in } (0,R).
%\ee
\begin{align}\label{ODEsystem}
u''+\frac{u'}{r}-\frac{k^2u}{r^2} &=\frac{u}{L}\left[-a^2+\sqrt{\frac{2}{3}} b^2 v+c^2\left( u^2+ v^2\right)\right],\nonumber\\
v''+\frac{v'}{r}&=\frac{v}{L}\left[-a^2-\frac{1}{\sqrt{6}}b^2 v+c^2\left( u^2+ v^2\right) \right] + \frac{1}{\sqrt{6} L} b^2 u^2, \ \ r \in (0,R). 
\end{align}
%Using boundary conditions at $R$, in order to obtain smooth solutions  and recalling that $Y$ has to solve \eqref{eq:EL}  we take the following boundary conditions on $(u,v)$
Boundary conditions at $r=0$ follow from requiring $Y$ to be a smooth solution of \eqref{eq:EL}, while boundary conditions at $r = R$ are given by \eqref{eq: bc at R}, as follows:
\be\label{bdrycond}
u(0)=0,\ v'(0)=0, \ u(R)=\frac{1}{\sqrt{2}} s_+,\,\,\,v(R)=-\frac{1}{\sqrt{6}}s_+.
\ee

In order to show  the existence of a  solution $Y$ of \eqref{eq:EL} of the form \eqref{ansatz}, 
we need to establish the existence of a solution of the system of ODEs \eqref{ODEsystem} -- \eqref{bdrycond}. We do this using methods of calculus of variations.
Substituting the ansatz \eqref{ansatz} into the Landau-de Gennes energy \eqref{FQ}, we obtain a reduced  1D energy functional corresponding to the system \eqref{ODEsystem}, 
\bea\label{def:mcR}
\mcE(u,v)=&\int_0^R \left[ \frac{1}{2} \left( (u')^2+(v')^2+\frac{k^2}{r^2}u^2 \right) -\frac{a^2}{2L}(u^2 + v^2)+\frac{c^2}{4L}\left(u^2+v^2\right)^2\right]\,rdr\non\\
&-\frac{b^2}{3L \sqrt{6} }\int_0^R v(v^2 -3u^2)\,rdr .
\eea
The energy $\mcE$ is defined on the admissible set 
\be\label{SR}
S = \left\{ (u,v) :  [0,R] \to \RR^2 \, \Big | \, \sqrt{r} u', \sqrt{r} v', \frac{u}{\sqrt{r}}, \sqrt{r} v \in L^2(0,R), \, u(R)=\frac{s_+}{\sqrt{2}}, v(R)= -\frac{s_+}{\sqrt{6}} \right\}.
\ee
\begin{theorem}\label{lemma:existenceODE}
For every $L>0$ and $0<R< \infty$, {  there exists a global minimiser $(u(r), v(r)) \in [C^\infty(0,R) \cap C([0,R])] \times [C^\infty(0,R) \cap C^1([0,R])] $ of the reduced energy \eqref{def:mcR} on $S$, which satisfies the system of ODEs \eqref{ODEsystem} -- \eqref{bdrycond}. } % {\color{red} which is a global minimiser of the reduced energy \eqref{def:mcR}}. 
\end{theorem}
\bproof
%Using  arguments from the calculus of variations we can show the existence of a minimiser of $\mcE$ on the set 
%\be\label{SR}
%S = \left\{ (u,v) :  [0,R] \to \RR^2 \, \Big | \, \sqrt{r} u', \sqrt{r} v', \frac{u}{\sqrt{r}}, \sqrt{r} v \in L^2(0,R), \, u(R)=\frac{s_+}{\sqrt{2}}, v(R)= -\frac{s_+}{\sqrt{6}} \right\}.
%\ee
%In fact 
It is straightforward to show that $\mcE(u,v) \geq -C$ for all $(u,v) \in S$. Therefore, there exists a minimizing sequence $(u_m, v_m)$ such that 
$$
\lim_{m \to \infty} \mcE(u_m,v_m) = \inf_S \mcE(u,v).
$$ 
Using the energy bound we obtain that $(u_m,v_m) \rightharpoonup (u,v)$ in $[H^1((0,R); r\,dr) \cap L^2((0,R); \frac{dr}{r})] \times H^1((0,R); r\,dr) $ (perhaps up to a subsequence). Using the Rellich-Kondrachov theorem and the weak lower semicontinuity of the Dirichlet energy term in $\mcE$, we obtain
$$
\liminf_{m \to \infty} \mcE(u_m,v_m) \geq  \mcE(u,v),
$$
which establishes the existence of a minimiser $(u,v) \in S$. 
{\color{black}
Since $(u,v)$ is a minimiser of $\mcE$ on $S$, it follows that $(u,v)$ satisfies the Euler-Lagrange equations  \eqref{ODEsystem}. Then the matrix-valued function  $Y:B_R(0)\to \mcS_0$  defined as in \eqref{ansatz} is a weak solution of the PDE system \eqref{eq:EL}, and thus is smooth and bounded on $B_R$ (see for instance \cite{mz}). Since $F_3$ is a constant matrix we have that  $v(r)=\tr(YF_3)\in C^\infty(0,R)\cap L^\infty(0,R)$. Similarily $F_n$ is smooth on $B_R\setminus\{0\}$ hence $u(r)=\tr(YF_n)\in C^\infty (0,R)\cap L^\infty (0,R)$.

Furthermore, since $u\in H^1((0,R); r\,dr) \cap L^2((0,R); \frac{dr}{r})$ we have for any $[a,b]\subset (0,R]$ that  $u\in H^1([a,b])$ hence continuous. Moreover, we have:
$$u^2(b)-u^2(a)=2\int_a^b u'(s)u(s)\,ds\le \left(\int_a^b |u'(s)|^2\,sds\right)^{\frac{1}{2}}\left(\int_a^b |u(s)|^2\,\frac{ds}{s}\right)^{\frac{1}{2}}.$$ 
Hence, taking into account that the right-hand side of the above tends to $0$ as $|b-a|\to 0$ we get that $u$ is continuous up to $0$ so $u\in C([0,R]) \cap L^2((0,R); \frac{dr}{r})$ and therefore $u(0)=0$.

Using the Euler-Lagrange equations for $v$ we obtain
$$
v'(r) =\frac{1}{r} \int_{0}^{r}  g ( {u} , {v} ) \, s \,ds , r>0
$$
where $g(u,v) = \frac{v}{L}\left[-a^2-\frac{1}{\sqrt{6}}b^2 v+c^2\left( u^2+ v^2\right) \right] + \frac{1}{ \sqrt{6}L} b^2 u^2  $. It follows that $\lim_{r \to 0} v'(r) =0$. Using again the equation for $v$  at $r=R$ we get that $v\in C^1([0,R])$.}
\eproof
\begin{remark}\label{remark:maxprinc}
Using maximum principle arguments  it is possible to show (see \cite{mz}) 
$$|Y|^2 = u^2+ v^2 \leq \frac{2}{3} s_+^2.$$
\end{remark}

%\bigskip
%We proved that there is a special solution $Y$ defined in \eqref{ansatz} satisfying Euler-Lagrange equations \eqref{eq:EL}. Now we would like to study properties of this solution when $f(Q)$ is Ginzburg-Landau type potential .
%
\section{The case $b=0$:  properties of $Y$}
\label{subsec:qualbzero}

In this section we concentrate on the problem \eqref{ODEsystem} for the case $b^2=0$.  In this case, the  bulk energy $f(Q)$ becomes the standard Ginzburg-Landau potential (that is, a double well potential in $|Q|^2$). We are then able to show that there is a unique global minimiser $(u,v)$ of the energy \eqref{def:mcR}, and that this minimiser satisfies $u>0$ and $v<0$ on $(0,R]$.
% (see Figure~\ref{fig:fig7}). 

%\begin{figure}[ht]
%\begin{minipage}[b]{0.45\linewidth}
%\centering
%\includegraphics[width=1\textwidth]{L1uk1k2.pdf}
%%\caption{default}
%\end{minipage}
%\hspace{0.5cm}
%\begin{minipage}[b]{0.45\linewidth}
%\centering
%\includegraphics[width=1.01\textwidth]{L1vk1k2.pdf}
%%\caption{default}
%\end{minipage}
%\caption{The profiles $u$ (left) and $v$  (right) for $R=L=1$, $k=1$ and $k=2$  }
%\label{fig:fig7}
%\end{figure}
\begin{lemma}\label{lemma:positivity}
  Let $L>0$, $0<R<\infty$, $b^2=0$.  Let $( u, v)$ be a global minimiser of \eqref{def:mcR} over the set $S$ defined in \eqref{SR}. 
Then:
  \begin{enumerate}
    \item  ${u} >0$ on $(0,R]$.
   {\color{black}  \item ${v} <0$ and $v' \ge 0$ on $[0,R]$. } % and ${v}' \geqslant 0$ on $[0,R]$.}
  \end{enumerate}
\end{lemma}

{\color{black}
{\bproof} We define $\tilde{u}  :=  | {u}
|$ and $\tilde{v} := - | {v} |$. We note that since $b^{2} =0$,  %the pair 
$(\tilde u, \tilde v)$ is a global minimiser of $\mcE$ on $S$.  It follows from Theorem~\ref{lemma:existenceODE} that
$\tilde u \in C^\infty(0,R) \cap C([0,R])$, $\tilde v \in C^\infty(0,R) \cap C^1([0,R])$ and that
$(\tilde u, \tilde v)$ satisfies the Euler-Lagrange equations  \eqref{ODEsystem} and boundary conditions \eqref{bdrycond} with $b^2 = 0$.

Suppose for contradiction that $\tilde u(r_0) = 0$ for some $r_0 \in (0,R)$.  Since $\tilde u$ is  smooth and nonnegative, it follows that $\tilde u'(r_0) = 0$.  On the other hand,  the unique solution of the  initial value problem for the second-order regular ODE  satisfied by $\tilde u$ (for given, fixed $\tilde v$): 
 \[ \tilde u''+\frac{\tilde u'}{r}-\frac{k^2\tilde u}{r^2} =\frac{\tilde u}{L}\left[-a^2+c^2\left( \tilde u^2+ \tilde v^2\right)\right]\]
on $(r_0,R)$ with initial conditions $u(r_0)=u'(r_0)=0$ is given by $\tilde u\equiv 0$ identically.  But this contradicts the fact that $\tilde u(R) = \frac{s_+}{\sqrt2}>0$.  Therefore, $\tilde u > 0$ on $(0,R)$, and since $u(R) > 0$, it follows that $u > 0$ on $(0,R]$.
 
 A similar argument shows that $v < 0$ on  $(0,R]$, which then allows us to establish that $v' \ge 0$ on $(0,R)$.  Indeed, from the Euler-Lagrange equation for $v$, it follows that
\[
  {v}' ( r ) = \frac{1}{r} \int_{0}^{r} \frac{v}{L}   \left[
   -a^{2} +c^{2}  ( u^{2} +v^{2} )\right ] 
  s  \,ds  .\]
From Remark~\ref{remark:maxprinc}, we get  that $u^2 + v^2 \le \frac{a^2}{c^2}$, which together with the preceding yields 
\[ v' \ge 0\text{  on } [0,R].\]
 Since $v(R) < 0$, it follows that $v(0) < 0$, so that $v<0$ on $[0,R]$.
 {\eproof}
 
% \begin{remark}
% From the Euler-Lagrange equation for $v$, we have that
%\[
%  {v}' ( r ) = \frac{1}{r} \int_{0}^{r} \frac{v}{L}   \left[
%   -a^{2} +c^{2}  ( u^{2} +v^{2} )\right ] 
%  s  \,ds  .\]
%  It follows from Remark~\ref{remark:maxprinc} and Lemma~\ref{lemma:positivity} that $v( -a^{2} +c^{2}  ( u^{2} +v^{2} ) )\ge 0$,
%  so that 
%  \[ v' \ge 0\text{  on } [0,R].\]
% \end{remark}
 }

%\[ \tilde u''+\frac{\tilde u'}{r}-\frac{k^2\tilde u}{r^2} &=\frac{\tilde u}{L}\left[-a^2+c^2\left( \tilde u^2+ \tilde v^2\right)\right]\]
%with initial conditions $\tilde u(r_0) = \tilde u'(r_0) = 0$ 
%sing standard 
%regularity arguments and the boundedness of $(u,v)$ (see Remark~\ref{remark:maxprinc}) we obtain that 
% $\tilde u$,$\tilde v\in C^2(0,R)$. The strong maximum principle applied to the equation for $\tilde u$ implies that $\tilde u(r)>0,\forall r\in (0,R)$. Since $u(R)>0$ we get $u(r)>0,\forall r\in (0,R)$. 
%Next we rearrange the second ODE in {\eqref{ODEsystem}} to obtain the equivalent form
%\begin{equation}
%  ( r_{} v' )' =r vf ( u,v ) , \hspace{1em} \text{with }  f ( u,v ) := \frac{1}{L}\left[
%  -a^{2} +c^{2}  ( u^{2} +v^{2} )\right ] \leqslant 0
%\end{equation} (see Remark~\ref{remark:maxprinc} and the definition of $s_+$ in \eqref{def:s_+}  for the last inequality).
%Integrating this expression on $[ 0,r ]$, we get that
%\begin{equation}
%  {v}' ( r ) = \frac{1}{r} \int_{0}^{r} {v}  f ( {u} , {v} ) 
%  s  \,ds    \geqslant   
% 0 .
%\end{equation}
%Since ${v}$ is a nondecreasing
% function and ${v} ( R )=-\frac{s_+}{\sqrt{6}} <0$,
%we get  that $v <0$ on $[ 0,R ]$. 
%}

%\newpage
\begin{proposition}\label{lemma:unique}
  Let $L>0$, $0<R<\infty$, $b^2=0$.  {There exists a unique {\color{black} solution of \eqref{ODEsystem}, \eqref{bdrycond} in the class of solutions satisfying $u>0,v<0$ on $(0,R)$.}} %\ Then the solution $( u, v)$ of \eqref{ODEsystem}, \eqref{bdrycond} obtained in Theorem~\ref{lemma:existenceODE} is unique in the class of solutions such that $u>0$ and $v<0$ on $(0,R]$. 
\end{proposition}
{\bproof}
{
Existence follows from Theorem~\ref{lemma:existenceODE} and Lemma~\ref{lemma:positivity}.  To prove uniqueness, we use the approach of Brezis and Oswald \cite{brezisoswald}. Suppose that $(u,v)$ and $(\xi, \eta) $ {\color{black} satisfy \eqref{ODEsystem}  with $u,\xi>0$ and $v,\eta<0$ on $(0,R)$.}
%different solutions of \eqref{ODEsystem}, \eqref{bdrycond}.  
We obtain
\bea
& \frac{\Delta_r u}{u} - \frac{\Delta_r \xi}{\xi} = \frac{1}{L} \left(  c^2(u^2 +v^2) - c^2 (\xi^2 + \eta^2) \right), \\
& \frac{\Delta_r v}{v} - \frac{\Delta_r \eta}{\eta} = \frac{1}{L} \left( c^2 (u^2 +v^2) - c^2 (\xi^2 + \eta^2) \right) ,
\eea
where $\Delta_r u = u'' + \frac{u'}{r}$.
%where $\Delta = d^2/dr^2 +1/r\, d/dr$.
Multiplying the first equation by $\xi^2 - u^2$ and the second equation by $\eta^2 -v^2$, and  then adding the two, we obtain
$$
\left( \frac{\Delta_r u}{u} - \frac{\Delta_r \xi}{\xi} \right) (\xi^2 - u^2) + \left(\frac{\Delta_r v}{v} - \frac{\Delta_r \eta}{\eta} \right) (\eta^2 -v^2) =
-\frac{c^2}{L}(u^2+v^2 - \xi^2 - \eta^2)^2.
$$
Multiplying by $r$, integrating over $[0,R]$ and taking into account that $u(R) =\xi(R)$,  $v(R) =\eta(R)$, we obtain
\begin{align*}
 & \int_0^R\left\{  \left[(u/\xi)'\xi\right]^2 + \left[(\xi/u)'u\right]^2 + \left[(v/\eta)'\eta\right]^2 + \left[(\eta/v)'v\right]^2 \right\} r \, dr \\
 &\qquad+ \int_0^R \frac{c^2}{L}(u^2+v^2 - \xi^2 - \eta^2)^2 \, r \, dr =0.
\end{align*}
%\bea
 %& \int_0^R \left\{ \left[u' - \frac{u}{\xi} \xi' \right]^2 +  \left[ \xi' - \frac{\xi}{u} u' \right]^2 +  \left( v' - \frac{v}{\eta} \eta' \right)^2 +  \left( \eta' - \frac{\eta}{v} v' \right)^2 \right\} r \, dr \\
 %&\qquad+ \int_0^R c^2(u^2+v^2 - \xi^2 - \eta^2)^2 \, r \, dr =0
%\eea
This implies $u(r) =k_1\xi (r) $ and $v(r) = k_2 \eta(r)$ for some $k_1,k_2\in\mathbb{R}$ and every $r \in [0,R]$. Therefore, due to the boundary conditions, we obtain $k_1=k_2=1$ and the proof is finished.
{\eproof}

%\begin{remark}
%Combining Lemma~\ref{lemma:positivity} and Proposition~\ref{lemma:unique} we deduce that for $b=0$ there is a unique global minimiser of the energy \eqref{def:mcR}.
%\end{remark}

Now we are ready to investigate the minimality of the %special 
solution of the Euler-Lagrange equation \eqref{eq:EL} introduced in section 3 % the previous sections 
with respect to variations $P \in H^1_0(B_R, \mcS_0)$. 
%In order to do this we study the sign of the second variation of the energy \eqref{FQ} at the point $Y$ defined in \eqref{ansatz}. It is straightforward to compute that  the second variation  is 
%\bea
%\CL_R[Y](P)&=\frac{1}{2}\frac{d^2}{dt^2}|_{t=0} \mathcal{\mcF}[Y+tP]\non\\
%&=\int_{B_R}\frac{1}{2}|\nabla P|^2-\frac{a^2}{2}|P|^2-b^2\tr(P^2Y)+\frac{c^2}{2}\left(|Y|^2|P|^2+2|\tr(YP)|^2\right),
%\eea 
%where $P \in H^1_0(B_R, \mcS_0)$.
We show that for $b^2=0$,  the solution $Y$ given by \eqref{ansatz} % of the Euler-Lagrange system \eqref{eq:EL}
 is the unique global minimiser of  energy \eqref{FQ}.
\begin{theorem} \label{th:stab}
Let $b^2=0$, and let $Y$ be given by \eqref{ansatz} with $(u,v)$ the unique global minimiser of the reduced energy \eqref{def:mcR} in the set  $S$ (defined in \eqref{SR}).
 Then $Y$ is the unique global minimiser of the Landau-de Gennes energy \eqref{FQ} in $H^1(B_R;\mcS_0)$. 
 \end{theorem}
 \bproof  
%In order 
%To prove the claim 
We take %$Y$ as in the conditions of the theorem,  
$P \in  H^1_0(B_R;\mcS_0)$ and compute the difference in energy between $Y+P$ and $Y$,  
% the energy difference %in  between two energies
\begin{equation}\label{eq: FPP}
  \begin{array}{lll}
    \mathcal{F} ( Y+P ) -\mathcal{F} ( Y ) & = & \mycal{\mathcal{I}} [ Y ] (
    P,P ) + \frac{1}{L}
     \int_{B_{R}} \frac{c^{2}}{4} ( | P |^{2} +2\,\tr(YP) )^{2} ,  
  \end{array}
\end{equation}
where
\begin{equation}\label{IPP}
  \mycal{\mathcal{I}} [ Y ] ( P,P ) =  \frac{1}{2} \int_{B_{R}} | \nabla
  P |^{2}  +   \frac{1}{2L}\int_{B_{R}} | P |^{2} \left( -a^2
 % \frac{a^{2}}{L} 
 + c^2 % \frac{c^{2}}{L}
  | Y |^{2} \right)   
\end{equation} and we have used the fact that $Y$ satisfies \eqref{eq:EL} in order to eliminate the first-order terms in $P$. 
%To conclude the proof 
Thus, it is sufficient to prove that
$\mycal{\mathcal{I}} [ Y ] ( P,P ) \geqslant C \| P \|_{L^2}$ for every $P \in H_{0}^{1}
\left( B_{R} ( 0 ) , \mycal{S_{0}} \right)$. % In order to prove this, 

%In order 
To investigate \eqref{IPP} we use a Hardy trick (see, for instance \cite{INSZ1}).
%We
%recall that $v \in C^{\infty} (0,R)$ satisfies the second
%equation of \eqref{ODEsystem}, and, 
From Lemma~\ref{lemma:positivity}, we have that $v<0$ on $[ 0,R ]$. 
Therefore, %It is clear that 
any $P \in H^1_0(B_R,\mcS_0)$ can be written in the form $P(x)= v(r) U(x)$, where $U \in  H^1_0(B_R,\mcS_0)$.  Using equation \eqref{ODEsystem} for $v$ we have the following identity
$$
  v \Delta v=
\frac{v^{2}}{L} \left( - {a^{2}} + {c^{2}} | Y |^{2} \right)
$$
and  therefore  
\bea
\mathcal{I} [ Y ] (  P,P )= \frac{1}{2}  \sum_{i,j} \int_{B_R} | \nabla v(|x|) U_{ij} (x) + v(|x|) \nabla U_{ij}(x) |^2 +  \Delta v(|x|) v(|x|) U_{ij}^2(x).
\eea 
Integrating by parts in the second term above, we obtain
$$
 \sum_{i,j} \int_{B_R} \Delta v \, v \, U_{ij}^2  = - \sum_{i,j} \int_{B_R}  |\nabla v|^2 U_{ij}^2 + 2 \nabla v \cdot \nabla U_{ij}\, v\, U_{ij} .
$$
It follows that
$$
\mathcal{I} [ Y ] ( P, P ) =  \frac{1}{2}  \int_{B_R} v^2 \, |\nabla U|^2 .
$$
Using the fact that $0< c_1 \leq v^2 \leq c_2$ (see Lemma~\ref{lemma:positivity}) and the Poincar\'e inequality we obtain
$$
\mathcal{I} [ Y ] ( P, P ) \geq C \int_{B_R} |P|^2.
$$
From \eqref{eq: FPP}, it follows that %Therefore we readily obtain 
\be\label{rel:coercivitybzero}
 \mathcal{F} ( Y+P ) -\mathcal{F} ( Y ) \ge C\|P\|^2_{L^2},
\ee %which, by standard arguments, suffices for obtaining the claimed stability.
therefore $Y$ is the 
%It is clear that $Y$ is the 
unique global minimizer of the energy $ \mathcal{F}$. % since $I[Y](P,P)>0$ for $P\not=0$. 
\eproof
\begin{remark}
It is  straightforward to use the continuity of the solutions $(u,v)$ with respect to the parameter $b^2$ to show that for $b^2$ small enough, the solution $(u_b,v_b)$ of \eqref{ODEsystem} -- \eqref{bdrycond} found in Theorem~\ref{lemma:existenceODE} generates a global minimizer $Y$ of the energy \eqref{FQ}. 
\end{remark}

\subsection{Limiting case $L \to 0$}
Next we  consider the limit $L \to 0$. We define the energy
$$
\mcE_L (u,v) = \int_0^R \left[  \frac{1}{2} \left( (u')^2+(v')^2+\frac{k^2}{r^2}u^2 \right) + 
\frac{c^2}{4L} \left( (u^2 + v^2) - \frac{a^2}{c^2} \right)^2 \right] r\, dr .
$$
For  $b=0$, $\mcE_L$ coincides with the reduced energy \eqref{def:mcR} up to an additive  constant. We also define the following space:
$$
H=\left\{ (u,v) :  [0,R] \to \RR^2 \, | \, \sqrt{r} u', \sqrt{r} v', \frac{u}{\sqrt{r}}, \sqrt{r} v \in L^2(0,R) \right\}.
$$
%We prove the following lemma
\begin{lemma}
In the limit $L \to 0$  the following statements hold:
\begin{enumerate}
\item If  $(u_L,v_L) \in S$ (see \eqref{SR}) and $\mcE_L (u_L,v_L) \leq C$, then $(u_L,v_L) \rightharpoonup (u,v)$  in $H$ (perhaps up to a subsequence).  Moreover, $(u,v) \in S$ and $u^2(r) + v^2(r) = \frac{a^2}{c^2}$ a.e. $r \in (0,R)$.
\item $\mcE_L $ $\Gamma$-converges to $\mcE_0 $  in $S$ , where
\be \label{eq: E_0}
\mcE_0 (u,v)= \left\{
\begin{array}{cl}
 \int_0^R \frac{1}{2} \left( (u')^2+(v')^2+\frac{k^2}{r^2}u^2 \right)  r\, dr &
\hbox{ if } u^2 + v^2 = \frac{a^2}{c^2}, \\
\infty  &  \hbox{ otherwise. }
\end{array} \right.
\ee
\end{enumerate}
\end{lemma}
\bproof
The first statement follows from the energy estimate $\mcE_L (u_L,v_L) \leq C$. % ({\color{red} see, for instance, \cite{mz}).}

Next we show the $\Gamma$-convergence result.  To do this we must  check the following:
\begin{itemize}
\item for any $(u_L, v_L) \in S$ such that {\color{black} $(u_L, v_L) \to (u,v)$ in $S$}, we have that
$$
\liminf_{L \to 0} \mcE_L (u_L,v_L) \geq \mcE_0 (u,v);
$$
\item for any $(u,v) \in S$, there exists a sequence  $(u_L, v_L) \in S$ such that
$$
\limsup_{L \to 0} \mcE_L (u_L,v_L) = \mcE_0 (u,v).
$$
\end{itemize}
The first part of  the $\Gamma$-convergence result follows from the lower semicontinuity of the Dirichlet term in the energy $\mcE_L$  and the penalization of the potential.
To prove the second part, we note that for any {\color{black} $(u,v) \in S$} we  may take the recovery sequence $(u_L,v_L) =(u,v)$, for which the  $\limsup$ equality is clearly satisfied. 
\eproof

Next we show that the global minimiser of $\mcE_0$ defines the unique global minimiser of a  certain harmonic map problem.% for $Q$. 
\begin{theorem} \label{th47}
Let $0<R<\infty$. There exist exactly two critical points of $\mcE_0$ over the set $S$ defined in \eqref{SR}. These are  explicitly given by the following formulae: 
\begin{align}\label{eq: u and v}
u_-(r) &= 2\sqrt{2} s_+ \, \frac{ R^{|k|} r^{|k|}} {r^{2|k|} + 3R^{2|k|}}, \quad v_-(r) = \sqrt{\frac23} s_+\,  \frac{r^{2|k|} - 3R^{2|k|}}{r^{2|k|} + 3R^{2|k|}},\nonumber\\
u_+(r) &= 2\sqrt{2} s_+ \, \frac{ R^{|k|} r^{|k|}} {3r^{2|k|} + R^{2|k|}}, \quad v_+(r) = \sqrt{\frac23} s_+\,  \frac{R^{2|k|} - 3r^{2|k|}}{3r^{2|k|} + R^{2|k|}}
\end{align} with $s_+$ given by \eqref{def:s_+} with $b^2=0$.
If we define 
$$
Y_{\pm} = u_{\pm} F_n + v_{\pm} F_3,
$$
then  $Y_-$ is the unique global minimiser and $Y_+$ is a critical point of the following harmonic map problem:
\be\label{hmp}
\min \left\{  \int_{B_R} \frac{1}{2} |\nabla Q|^2 \, \Big | \, Q \in H^1 (B_R, \mcS_0), \, Q(R) = Q_k ,\, |Q|^2=\frac{2}{3}s_+^2 \textrm{ a.e.  in  }B_R \right\}.
\ee
\end{theorem}
\bproof
{
The constraint $u^2 + v^2 = \frac{a^2}{c^2}$ may be incorporated through the substitution
\begin{equation}
\label{eq: psi }
u = \sqrt{\frac23} s_+ \sin\psi, \ \ v = -\sqrt{\frac23}s_+ \cos\psi,
\end{equation}
where $\psi: (0,R] \rightarrow \RR$.
In terms of $\psi$, the energy $\mcE_0$ is given up to a multiplicative constant by
\begin{equation}
\label{eq: E_0 psi }
\mcE_0[\psi] =  \frac12 \int_0^R \left(r {\psi'}^2 + \frac{k^2}{r}\sin^2\psi \right) \,dr.
\end{equation}
Critical points of $\mcE_0$ satisfy the Euler-Lagrange equation
\begin{equation}
\label{eq: EL 2 }
\left(r\psi'\right)' = \frac{k^2}{r} \sin \psi \cos \psi,
\end{equation}
and therefore belong to $C^\infty(0,R)$.  From \eqref{eq: bc at R} and \eqref{eq: psi }, $\psi$ satisfies
 the boundary condition $\psi(R) = \frac{\pi}{3}+ 2 \pi j$ for  $ j \in \ZZ$. Without loss of generality, we may take $j = 0$ (since $\psi$ and $\psi + 2\pi j$ correspond to the same $(u,v)$).  Therefore, we may take the boundary condition as
 \begin{equation}
\label{eq: psi bc}
\psi(R) = \frac{\pi}{3}.
\end{equation}

The Euler-Lagrange equation \eqref{eq: EL 2 } may be integrated to obtain the relation
\be \label{eq: first integral} 
\frac12 r^2 {\psi'}^2 - \frac{k^2}{2}  \sin^2 \psi = -\frac{k^2}{2} \alpha \ee
for some constant $\alpha \le 1$.  We claim that $\alpha = 0$.  
First, we note that  $\alpha < 0$ would imply that    $r^2 {\psi'}^2$ is bounded away from zero, which is incompatible with $\mcE_0[\psi]$ being finite.  Next, $\alpha = 1$ would imply that
 $\sin^2 \psi = 1$ identically, which is incompatible with the boundary condition \eqref{eq: psi bc}.  It follows that $0 \le \alpha < 1$.
If $\alpha > 0$, we may define $x(t) = \psi(\exp t)$ for $t \in (-\infty, \ln R)$.  Then
$\frac12 {\dot x}^2 = \frac{k^2}{2}\left(\sin^2 x - \alpha\right)$.
It is an elementary result (the simple pendulum problem) that $x(t)$ is periodic with period $T$ (we omit the explicit expression for $T$); this implies that  $\psi \left( e^{-T} r\right)  = \psi(r)$.  In addition, $A:= \int_{\tau}^{\tau + T} \sin^2 x \, dt $ is strictly positive and independent of $\tau$; in terms of $\psi$,  this implies that
\[ \int_{e^{-nT} R}^R \frac{\sin^2 \psi}{r}\, dr = nA\]
for $n \in \mathbb{N}$.  It follows that $u^2/r = \frac23 s_+^2 \sin^2 \psi/r $ is not square-integrable, which is incompatible with $\mcE_0[\psi]$ being finite.  Thus we may conclude that $\alpha = 0$.  

We claim now that any solution of  \eqref{eq: EL 2 }   satisfies either
$r \psi'(r) =  |k| \sin \psi$ or $r \psi'(r) =  -|k| \sin \psi$ on the whole interval $(0,R)$. For suppose  $\chi (r)$ is a smooth solution of \eqref{eq: EL 2 }, and that for some point $r_0\in (0,R)$ we have that $r_0 \chi '(r_0) = |k| \sin \chi(r_0)$. Then regarding \eqref{eq: EL 2 } as a {\it regular} second-order ODE on $(0,R)$ we have that the {\it initial-value problem} \eqref{eq: EL 2 } with initial conditions $\psi(r_0)=\chi(r_0)$, $\psi'(r_0)=\frac{|k|}{r_0}\sin \chi(r_0)$ has a unique smooth solution on $(0,R)$, namely the one satisfying the first order equation $\chi'(r)=\frac{|k|}{r}\sin\chi(r)$ on $(0,R)$, which proves our claim.

Solving the first-order separable ODEs and applying the boundary conditions \eqref{eq: psi bc} we obtain exactly two solutions $\psi_\pm$ satisfying
$$
\tan \frac{\psi_\pm(r)}{2} =\frac{1}{ \sqrt{3}} \left(\frac{r}{R}\right)^{\mp |k|}.
$$
These correspond via \eqref{eq: psi } to \eqref{eq: u and v}.}

It is straightforward to check using the definition of $Y_\pm$ and \eqref{deltatopolar} that 
$$
\Delta Y_\pm = -\frac{3}{2s_+^2}|\nabla Y|^2Y_\pm, \ |Y_\pm|^2 =\frac{2}{3} s_+^2, \ Y_\pm(R, \varphi) =Q_k (\varphi).
$$
Therefore, $Y_\pm$ are critical points of the harmonic map problem \eqref{hmp}. 

Next, we show that  $Y_-$ is the unique global minimiser of the harmonic map problem \eqref{hmp}. 
Take $P \in H^1_0(B_R; \mcS_0)$ such that $|Y_-+P|^2=\frac23 s_+^2$. Then
$$
\frac12 \int_{B_R} |\nabla (Y_-+P)|^2 - \frac12 \int_{B_R} |\nabla Y_-|^2 = \frac12 \int_{B_R} |\nabla P|^2 + 2 \sum_{ij} \nabla [Y_{-}]_{ij} \cdot \nabla P_{ij}.
$$
Integrating by parts and using the Euler-Lagrange equation for $Y_-$, we obtain
$$
\frac23 s_+^2 \int_{B_R} \sum_{ij} \nabla [Y_{-}]_{ij}  \cdot \nabla P_{ij} =  \int_{B_R}  |\nabla Y_-|^2 \, \tr(Y_-  P).
$$
Using the fact that $|P|^2 = -2\, \tr (Y_-  P)$ we obtain
$$
\frac12 \int_{B_R} |\nabla (Y_-+P)|^2 - \frac12 \int_{B_R} |\nabla Y_-|^2 =  \frac12 \int_{B_R} |\nabla P|^2  - \frac{3}{2s_+^2}   |\nabla Y_-|^2 |P|^2 .
$$
The fact that $Y_-$ is harmonic implies that 
$$
\Delta v_- = -\frac{3}{2s_+^2}\,  v_-  |\nabla Y_-|^2,
$$
and we have that $v_- <0$  on $[0,R]$. Therefore
$$
\frac12 \int_{B_R} |\nabla (Y_-+P)|^2 - \frac12 \int_{B_R} |\nabla Y_-|^2 =  \frac12 \int_{B_R} |\nabla P|^2  + \frac{\Delta v_-}{v_-} |P|^2
$$
Using the decomposition $P = v(r) U$ and applying the Hardy decomposition trick in exactly the same way as in the proof of Theorem~\ref{th:stab}, we obtain
 $$
 \frac12 \int_{B_R} |\nabla (Y_-+P)|^2 - \frac12 \int_{B_R} |\nabla Y_-|^2 \geq C \| P\|^2_{L^2}
 $$
Therefore $Y_-$ is unique global minimiser of harmonic map problem \eqref{hmp}.
\eproof
\begin{remark}
It is straightforward to check that in the limit  $L \to 0$, the $\Gamma$-limit of the Landau-de Gennes energy
$$
\mcF(Q) = \frac12 \int_{B_R} |\nabla Q|^2 + \frac{c^2}{4L} \left( |Q|^2 - \frac23 s_+^2 \right)^2
$$
is exactly the harmonic map problem \eqref{hmp}.
\end{remark}
\begin{remark} \label{rm49} For $k$ even,
there is another explicit solution of the harmonic map problem \eqref{hmp}.  Let
%We restrict ourselves to {\it even $k$} and uniaxial $Q$-tensors 
\be \label{eq: uniaxial Q}
U=s_+ \left( m \otimes m - \frac{1}{3} I \right),
\ee
where
$$
m(r,\phi) =  \left( \frac{2 R^{\frac{k}{2}}  r^{\frac{k}{2}}}{ R^k + r^k} \cos\left(\frac{k\phi }{2}\right), \frac{2 R^{\frac{k}{2}}  r^{\frac{k}{2}}}{ R^k + r^k} \sin\left(\frac{k\phi }{2}\right), \frac{R^k - r^k}{ R^k + r^k} \right) .
$$
We note that $U$ is {\it uniaxial} (i.e., two of its eigenvalues are equal).
It is straightforward to show that %uniaxial ansatz 
$U$ is a critical point of the harmonic map problem \eqref{hmp}. Computing energies of $Y_+$, $Y_-$ and $U$ explicitly, we obtain
$$\color{black}
\mcE_D (Y_-) = \frac{2}{3} |k| \pi s_+^2 < 2 |k| \pi s_+^2 = \mcE_D(Y_+) =\mcE_D(U),
$$
where $\mcE_D (Q) =\frac12 \int_{B_R} |\nabla Q|^2$ is the Dirichlet energy.
\end{remark}

\begin{remark}\label{rm410}
The harmonic map \eqref{eq: uniaxial Q} % in the preceding remark 
is an example of a more general construction.  Let $\zeta = x + iy$, and let $f(\zeta)$ be {\color{black} meromorphic}.  Let
$$
{\color{black} m(x,y)} = \frac{ \left(2 \text{Re}\, f, 2 \text{Im}\, f, 1 - |f|^2\right)}{1 + |f|^2}.
$$
Then it is straightforward to show that $m$ defines an $S^2$-valued harmonic map (note that $|m| = 1$), and that ${\color{black} U := \sqrt{3/2} (m\otimes m - \frac13 I)}$ defines an $S^4$-valued harmonic map.  The example \eqref{eq: uniaxial Q}  is obtained by taking ${\color{black} f = (\zeta/R)^{k/2}}$, which corresponds to the boundary conditions \eqref{BC}.
\end{remark}

\begin{remark} \label{m411}
\color{black}
The results of \cite{danpatty} imply that for $|k| > 1$  and $b^2>0$ the global minimiser $Y$ of a reduced energy    in the limit $L \rightarrow 0$ approaches a harmonic map different  from $Y_-$.  In that case, the limiting harmonic map has $|k|$ isolated defects of index ${\mathop{\rm sgn}(k)}/2$.  \end{remark}

\section{Conclusions and outlook}
We have found a new highly symmetric equilibrium solution $Y$ of the Landau-de Gennes model, relevant for the study of liquid crystal defects  %defined in
of the form \eqref{ansatz}. This solution is valid for all values of parameters $a,b,c$, elastic constant $L$ and index $k$. The properties of this solution can be explored by investigating the system of ordinary differential equations \eqref{ODEsystem} -- \eqref{bdrycond}. 

We have provided a detailed study of solution $Y$ in the deep nematic regime when the material parameter $b^2$ is small enough (see \cite{coaxialvirga, MakGar} for a discussion on the physical relevance of this regime). In this case we have shown that  $Y$ is a global minimiser of the Landau-de Gennes energy, provided $(u,v)$ is a global minimiser of the energy \eqref{def:mcR}. In this sense, we have constructed the unique ground state of the 2D point defect, and linked its study to analyzing solutions of the ordinary differential equations \eqref{ODEsystem} --\eqref{bdrycond}. 

In the limiting case $L \to 0$ {\color{black} for $b^2 = 0$}, we have obtained for all $k$ two explicit defect profiles $Y_-$ and $Y_+$ (see Figure~\ref{fig:fig7}),  defined in Theorem~\ref{th47}.  The global minimiser $Y$ is equal to $Y_-$.  For even $k$, we obtain a third explicit  profile $U$ (see Figure~\ref{fig:fig9})
%-- $Y_-$, $Y_+$ and $U$, defined in Theorem~\ref{th47} and
defined in Remark~\ref{rm49}. % , respectively. 
It is straightforward to compute the eigenvalues of $Y_\pm$ and $U$ (see Figure~\ref{fig:fig8}),
\bea
\lambda^\pm_1 &= \sqrt{\frac{2}{3}} v^\pm(r), \  \ \ \lambda^\pm_2= -\frac{u^\pm}{\sqrt{2}} - \frac{v^\pm}{\sqrt{6}}, \ \ \ \lambda^\pm_3 = \frac{u^\pm}{\sqrt{2}} - \frac{v^\pm}{\sqrt{6}}, \\
\lambda^U_1 &= \lambda^U_2 =-\frac13, \  \ \ \lambda^U_3=\frac23.
\eea

\begin{figure}[ht]
\begin{minipage}[b]{0.45\linewidth}
\centering

%%%%%REMOVE COMMENT
\includegraphics[width=1\textwidth]{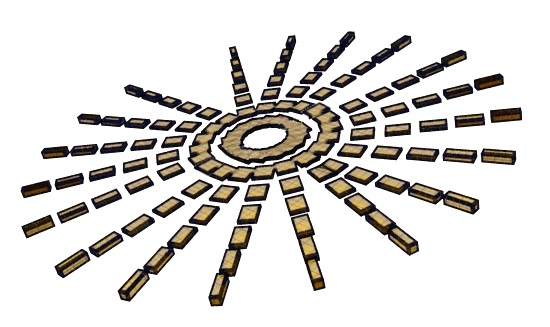}
%\caption{default}
\end{minipage}
\hspace{0.5cm}
\begin{minipage}[b]{0.45\linewidth}
\centering
%%%%%REMOVE COMMENT
\includegraphics[width=1\textwidth]{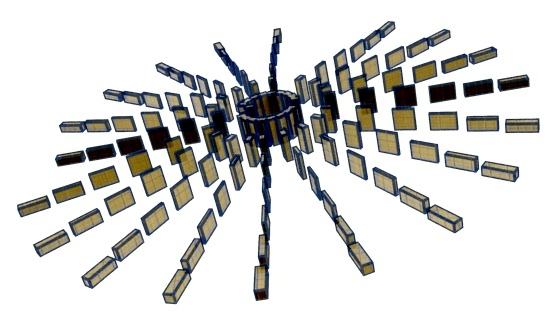}
%\caption{default}
\end{minipage}
\caption{$Y_-$ (left) and $Y_+$ (right) defects for of strength $1$   }
\label{fig:fig7}
\end{figure}

\begin{figure}[ht]
\begin{minipage}[b]{0.45\linewidth}
\centering
%%%%%REMOVE COMMENT
\includegraphics[width=0.9\textwidth]{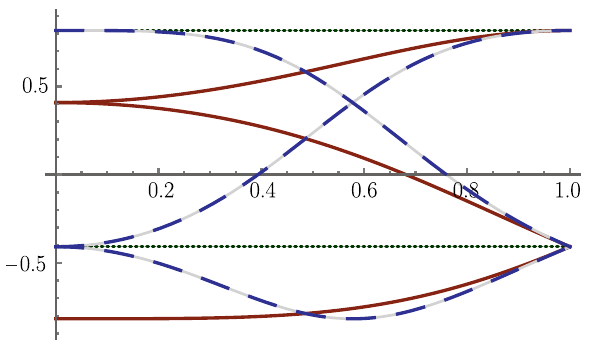}
%\caption{default}
\caption{Eigenvalues of $1$-strength defects: $Y_-$ (solid), $Y_+$ (dashed), $U$ (dotted) }
\label{fig:fig8}
\end{minipage}
\hspace{0.5cm}
\begin{minipage}[b]{0.45\linewidth}
\centering
%%%%%REMOVE COMMENT
\includegraphics[width=1\textwidth]{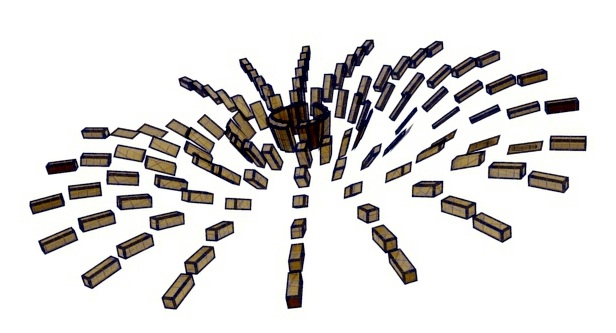}
%\caption{default}
\caption{Uniaxial defect of strength $1$ }
\label{fig:fig9}
\end{minipage}
\end{figure}

%Let's discuss the case $k=1$ in more detail. 
% The discussion isn't specific to k = 1, and the figures are for k = 1 and k = 2.
It is clear that the global minimiser $Y_-(r)$ is always biaxial except for points $r=0$ and $r=R$, while the critical point $Y_+$ is uniaxial at $0$, $R$ and the point of intersection of $\lambda_1^+$ and $\lambda_2^+$. Moreover, it is clear that $\lambda_3$ is the smallest eigenvalue. The structure of the defect profile $Y_+$ bears a resemblance to the three-dimensional {\it biaxial torus} profile \cite{MakGar}. However, whereas the biaxial torus is a candidate for the ground state in three dimensions, in this two-dimensional setting $Y_+$ has higher energy than $Y_-$, at least  in the small-$L$ regime. The profile $U$ is always uniaxial and  its energy coincides with the energy of $Y_+$.  

It is a very interesting and challenging task to find the ground state and  universal profile of the 2D defect for general parameters $a,b,c, L$. We are planning to tackle this problem  in the future.

\section*{Acknowledgement}
GDF, JMR, VS would like to acknowledge support from EPSRC grant EP/K02390X/1. VS also acknowledges support from EPSRC grant EP/I028714/1. AZ gratefully acknowledges the hospitality of the Mathematics Department at the University of Bristol, through EPSRC grants EP/I028714/1 and EP/K02390X/1.

\end{document}